\documentclass[12pt]{article}
\usepackage{amsmath, amssymb, amsthm,epsfig, latexsym }
\begin{document}
\begin{centering}

{\bf Some new behaviour in}\\ {\bf the deformation theory of Kleinian groups}

\vspace{3mm}

John Holt \\

\vspace{5mm}
\end{centering}
\begin{abstract}
Anderson and Canary (\cite{Anderson/Canary}) have constructed examples of 3-manifolds $M$ for which the space of convex co-compact representations of $\pi _1(M)$ in $PSL_2(\mathbb{C})$ is disconnected but has connected closure (in the topology of algebraic convergence).  More precisely, they showed that for their examples the intersection of the closures of any two path components of the space of convex co-compact representations is non-empty.  We study these examples and show that there is a connected, uncountable set of geometrically finite representations which is contained in the closure of every path component of the space of convex co-compact representations.
\end{abstract}
\vspace{3mm}

{\bf Introduction}

\vspace{2mm}

Suppose $M$ is an irreducible, orientable, atoroidal, compact 3-manifold whose boundary is non-empty.  Thurston's Hyperbolization Theorem (\cite{Morgan}) guarantees that the interior of $M$ may be imbued with the structure of a geometrically finite hyperbolic manifold.  Moreover, if $\partial M$ contains no torus components this structure may be taken to be convex co-compact.  But in particular there is a discrete, faithful representation of the fundamental group of $M$ into $PSL_2(\mathbb{C})$, the group of orientation preserving isometries of hyperbolic 3-space, $\mathbb{H}^3$, modelled as the upper half space in $\mathbb{R}^3$.  If $\Gamma$ is the image of this representation (so that $int(M)\cong \mathbb{H}^3/\Gamma$) then conjugating $\Gamma$ by an element of $PSL_2(\mathbb{C})$ produces a group whose quotient manifold is isometric to $\mathbb{H}^3/\Gamma$.  Hence if we are to understand all the possible hyperbolic structures on the interior of $M$ up to isometry we should study the set $H(\pi_1(M))$ of equivalence classes of faithful representations of $\pi _1(M)$ in $PSL_2 (\mathbb{C})$ whose images are discrete, where two representations are in the same equivalence class if they are conjugate in $PSL_2(\mathbb {C})$.

\vspace{2mm}

$H(\pi _1 (M))$ is a quotient of $D(\pi _1(M))=\{\rho \in Hom(\pi _1(M),PSL_2 (\mathbb{C}))|\rho $ is faithful and has discrete image$\}$ and so we can topologize $H(\pi _1 (M))$ by giving $D(\pi _1(M))$ the compact-open topology and $H(\pi _1(M))$ the quotient topology.  This is called the algebraic topology on $H(\pi _1(M))$.  Call $H(\pi_1(M))$ with the algebraic topology $AH(\pi_1 (M))$.  An element of $AH(\pi _1(M))$ does not necessarily give a hyperbolic structure to the interior of $M$; it may give a hyperbolic structure to a compact manifold homotopy equivalent to $M$ or the quotient of its image may not admit a manifold compactification, though it is a conjecture of Marden's that the latter does not occur.  Let ${\cal A}(M)$ denote the set of marked homeomorphism classes of compact, oriented, atoroidal, irreducible 3-manifolds  homotopy equivalent to $M$.  It follows from work of Ahlfors, Bers, Kra, Marden, Maskit, Sullivan and Thurston that the components of the interior of $AH(\pi _1(M))$ are enumerated by the elements of ${\cal A}(M)$ (see \cite{Canary/McCullough}).  In particular the homeomorphism type is constant on each component of the interior.  In the case when $\pi _1(M)$ is the fundamental group of a surface Bers (\cite{Bers}) conjectured that any B-group is contained in the closure of some Bers slice, and later Thurston (\cite{Thurston}) and Sullivan (\cite{Sullivan}) extended this to conjecture that the interior of $AH(\pi _1(M))$ is dense in $AH(\pi _1(M))$.

\vspace{2mm}

We specialize to the case when $\partial M$ contains no tori (so that the interior of $AH(\pi _1(M))$ is $CC(\pi _1(M))$, the subspace consisting of all those representations whose image is convex co-compact, see \cite{Marden} and \cite{Sullivan}).  Anderson and Canary (\cite{Anderson/Canary}) have constructed a family of examples for which there are finitely many components of $CC(\pi _1(M))$ and the closures of any two such components intersect.  They exhibit, for each integer $k\geq 3$, a manifold $M_k$ such that ${\cal A}(M_k)$ has $(k-1)!$ elements and they construct for each $[(M,h)] \in {\cal A}(M_k)$ representations $\rho _i$ in the component of $CC(\pi _1(M_k))$ indexed by $[(M_k, id)]\in {\cal A}(M_k)$ so that $\rho _i$ converges to $\rho$ and $\rho$ is in the closure of the component of $CC(\pi _1(M_k))$ indexed by $[(M,h)]$ (in fact $\rho$ uniformizes $M$).  This is to say, the homeomorphism type can change in the limit and the closure of $CC(\pi _1(M_k))$ is connected. (Recently Anderson, Canary and McCullough (\cite{Anderson/Canary/McCullough}) have generalized this result to give necessary and sufficient conditions for the closures of two components of the interior of $AH(\pi _1(M))$ to intersect, where now $M$ is any compact, orientable, irreducible, atoroidal 3-manifold with incompressible boundary.)

\vspace{2mm}

In this paper we examine Anderson and Canary's examples and show that for these examples there is a set of geometrically finite representations which is contained in the closure of every component of $CC(\pi _1(M_k))$.  The construction shows that in some sense this set can be chosen to be ``large'', in that given any $K\geq 1$ the set contains all of the $K$-quasiconformal deformations of some geometrically finite representation of $\pi _1(M)$ (this representation depending on $K$).

\vspace{2mm}

In the first section we will describe the construction of the examples and state the theorems that will be proven in section III. In section II we will present the versions of the Klein-Maskit combination theorems that we will be using, a statement of a generalization of the Hyperbolic Dehn Surgery theorem due to Comar (upon which the proofs in section III are largely based) and a few technical lemmas that we will use in section III.   

\vspace{2mm}
{\bf Acknowledgements}	This paper represents a portion of the author's Ph.D. thesis.  The author would like to thank his advisor, Dick Canary, for his advice and encouragement, as well as Richard Evans for many enjoyable discussions.

\vspace{3mm}

{\bf I:  The examples}

\vspace{2mm}

Choose an integer $k\geq 3$.  Let  $F(j)$ be a surface of genus $j$ with one boundary component and set $B(j)=F(j)\times [0,1]$, $j=1,...,k$.  Let $V=D^2 \times S^1$ be a solid torus and let $A(j)$ ($j=1,..,k$) be the annulus on $\partial V$ given by $A(j)=[e^{2\pi i (4j-1)/4k}, e^{2\pi i (4j+1)/4k}]\times S^1$.  Set $\partial _0 B(j)=\partial F(j)\times [0,1]$.  Form $M_k$ by identifying $\partial _0 B(j)$ with $A(j)$, $j=1,..,k$, with orientation-reversing homeomorphisms.  In other words, $M_k$ is obtained by attaching, in order, the $I$-bundles $B(1),\ldots,B(k)$ to a solid torus $V$ by identifying the annuli $\partial _0 (B(1)),\ldots, \partial _0 (B(k))$ with a collection of $k$ disjoint, parallel, longitudinal annuli $A(1),\ldots, A(k)\subset \partial V$. 

\vspace{2mm}

Suppose $\tau \in S_k$ is a permutation of the integers $1,...,k$.  Then we may obtain a manifold homotopy equivalent to $M_k$ by performing the above sort of construction but sewing $\partial _0 B(\tau(j))$ to $A(j)$.  Denote the manifold so obtained by $M_k ^{\tau}$.  $M_k^{\tau}$ is homeomorphic to $M_k^{\tau '}$ if an only if $\tau$ and $\tau '$ belong to the same right coset of $D_k$, the dihedral subgroup of $2k$ elements of $S_k$.  

\vspace{2mm}

Let $C_1 = \{e^{3 \pi i / 4k}  \}\times S^1 \subset \partial V$ and $C_2=\{e^{5\pi i /4k} \}\times S^1 \subset \partial V$ be two parallel curves in $\partial V\cap \partial M_k ^{\tau}$.  Form $\widehat {M_k ^{\tau}}$ by attaching $S^1\times [0,1]\times [0,1]$ to $M_k ^{\tau}$ by an embedding $h:S^1\times [0,1]\times \{0,1\}\longrightarrow \partial V$ such that $h(S^1 \times \{ \frac{1}{2}\}\times \{0\})=C_1$ and $h(S^1 \times \{ \frac{1}{2}\}\times \{1\})=C_2$.  $\widehat {M_k ^{\tau}}$ is homeomorphic to the manifold obtained by removing an open neighbourhood of the core curve of $V$ from $M_k ^{\tau}$.
 
\vspace{2mm}

Let $h_{\tau}:M_k \longrightarrow M_k ^{\tau}$ be a fixed homotopy equivalence that is the identity when restricted to the solid torus $V$.  Anderson and Canary have proven that $\{(M_k ^{\tau}, h_{\tau})| \tau \in S_k/\mathbb{Z}_k\}$ is a complete set of representatives for ${\cal A}(M_k)$.  For our purposes it is more convenient (and more natural in light of Johannson's Deformation Theorem(\cite{Johannson})) to require that $h_{\tau}$ be a homotopy equivalence between $M_k$ and $M_k ^{\tau}$ which is the identity {\underline {off}} of the solid torus $V$.  Using essentially the same proof as in ~\cite{Anderson/Canary} we have the following:

\vspace{2mm}

{\bf Lemma}(3.2 in \cite{Anderson/Canary}). {\it Let $\{\tau_1,...,\tau_{N}\}$ be a set of right coset representatives of $\mathbb{Z}_k = <(123\cdots k)>$ in $S_k$ ($N=(k-1)!$).  Then $\{(M_k ^{\tau _1},h_{\tau_1}),...,(M_k ^{\tau_N},h_{\tau_N})\}$ is a complete indexing set for the set of path components of $CC(\pi _1(M_k))$.}

\vspace{3mm}

In this paper we prove the following two results:

\vspace{3mm}

{\bf Theorem A  }{\it  For any integer $k\geq 3$ there is a representation which lies in the closure of every component of $CC(\pi _1(M_k))$.}

\vspace{3mm}

{\bf Theorem B  }{\it For any integer $k\geq 3$ and  $K\geq 1$ there is a representation $\rho _K$ such that the set of representations of $\pi _1 (M_k)$ induced by $K$-quasiconformal deformations of $\rho _K$ is contained in the closure of every component of $CC(\pi_1(M_k))$.}

\vspace{5mm}

{\bf II: Preliminaries}

\vspace{3mm}

{\bf Definitions}

\vspace{2mm}

In this article a  {\it Kleinian group} is a discrete, torsion-free subgroup of the group of isometries of hyperbolic 3-space, $\mathbb{H}^3$.  Realizing hyperbolic space as the upper half space $\{ (z,t) \ | \ t>0, \ z \in \mathbb{C}\}$ allows us to view a Kleinian group as a subgroup of $PSL_2(\mathbb{C})$.  If $\Gamma$ is a Kleinian group then the orbit space $\mathbb{H}^3/\Gamma$ is an orientable 3-manifold with constant sectional curvature equal to $-1$.  By fixing an orientation on ${\mathbb H}^3$ and passing this orientation to the quotient we may assume that the orbit space of a Kleinian group is an oriented 3-manifold.

\vspace{2mm}

If $M$ is a compact, oriented 3-manifold and $\Gamma$ is a Kleinian group such that $\mathbb{H}^3/\Gamma$ is homeomorphic to $int(M)$ via an orientation preserving homeomorphism we say that $M$ is {\it  uniformized } by $\Gamma$.  If $\Gamma$ is also the image of a representation $\rho:\pi _1(M) \longrightarrow PSL_2(\mathbb{C})$ then we sometimes say that $\rho$ uniformizes $M$.

\vspace{2mm}

The {\it limit set} of a Kleinian group $\Gamma$ is the set of accumulation points in $\widehat {\mathbb C}$ for the action of $\Gamma$ on $\mathbb {H}^3$.  The limit set of $\Gamma$ is denoted by $\Lambda (\Gamma)$.  When $\Gamma$ is non-abelian $\Lambda (\Gamma)$ is the smallest, non-empty, closed, $\Gamma-$invariant subspace of $\widehat {\mathbb C}$.  The complement of $\Lambda (\Gamma)$ in $\widehat {\mathbb{C}}$ is the {\it domain of discontinuity} for $\Gamma$ and is denoted by $\Omega (\Gamma)$.  As is implied by the name $\Gamma$ acts properly discontinuously on $\Omega (\Gamma)$.

\vspace{2mm}

A {\it Fuchsian group} is a Kleinian group whose limit set is a circle  and the stabilizer of any component of its domain of discontinuity is the whole group.  (Sometimes such a group is called Fuchsian of the {\it first kind}.)  For our purposes a Fuchsian group has the extended real line as its limit set.  That is, for our purposes a Fuchsian group is a subgroup of $PSL_2(\mathbb R)$.
\vspace{2mm}

The {\it convex core} of a hyperbolic 3-manifold $N= \mathbb{H}^3/\Gamma$ is the smallest convex sub-manifold of $N$ whose inclusion map is a homotopy equivalence.  The convex core of $N$ is denoted by $C(N)$.  Equivalently, $C(N)$ is the quotient of the convex hull of $\Lambda (\Gamma)$ in $\mathbb{H}^3$ by $\Gamma$.  If $C(N)$ is compact, then we say that $N$ (or $\Gamma$) is {\it convex co-compact}.  If $C(N)$ has finite volume and $\pi _1(N)$ is finitely generated then we say that $N$ is {\it geometrically finite}.  

\vspace{2mm}

For a hyperbolic 3-manifold $N=\mathbb{H}^3/\Gamma$ the {\it conformal extension} of $N$ is the manifold with boundary given by $(\mathbb{H}^3 \cup \Omega (\Gamma))/\Gamma$. 

\vspace{2mm}

Given a Kleinian group $\Gamma$ and a subgroup $J\subset \Gamma$ a subset $B\subset \widehat {\mathbb{C}}$ is {\it precisely $J$-invariant in $\Gamma$} (or, equivalently {\it precisely invariant under $J$ in $\Gamma$}), if $B$ is invariant under $J$ and if $g\in \Gamma$ is such that $g(B)\cap B\neq \emptyset$ then $g\in J$.

\vspace{2mm} 

A {\it quasiconformal deformation} of a representation $\rho$ is a representation $\rho '$ such that there is a quasiconformal map $f$ of $\widehat {\mathbb{C}}$ with $\rho '= f\circ \rho \circ f^{-1}$.  A {\it quasifuchsian group} is a quasiconformal deformation of a Fuchsian group.

\vspace{2mm}

For a compact, atoroidal, irreducible, orientable 3-manifold $M$ we have the set ${\cal A}(M)$ of homeomorphism classes of marked, compact, atoroidal, oriented, irreducible 3-manifolds homotopy equivalent to $M$.  Explicitly ${\cal A}(M)$ consists of equivalence classes of pairs $[(M',h)]$, where $h:M \longrightarrow M'$ is a homotopy equivalence and two pairs $(M_0, h_0)$ and $(M_1,h_1)$ are in the same equivalence class if and only if there exists an orientation preserving homeomorphism $j:M_0\longrightarrow M_1$ such that $j\circ h_0$ is homotopic to $h_1$.

\vspace{2mm}
$M_0$ is a {\it compact core} for an irreducible 3-manifold $M$ if $M_0$ is a compact, co-dimension 0 submanifold of $M$ whose inclusion map is a homotopy equivalence.  By Scott \cite{Scott} any 3-manifold with finitely generated fundamental group has a compact core.  It is a theorem of McCullough, Miller and Swarup (\cite{McCullough/Miller/Swarup}) that if $i_1:M_1 \rightarrow M$ and $i_2:M_2 \rightarrow M$ are two compact cores for an oriented 3-manifold $M$ then $i_1\circ \overline{i_2}$ is homotopic to an orientation preserving homeomorphism, where $\overline{i_2}$ is a homotopy inverse for $i_2$.

\vspace{2mm}

We have a map $\Theta:AH(\pi_1 (M))\longrightarrow {\cal A } (M)$ defined as follows.  Let $\rho \in AH(\pi _1(M))$ and let $M_{\rho}$ be a compact core for ${\mathbb H}^3/\rho (\pi _1(M))$.  Let $i_{\rho}$ be inclusion of $M_{\rho}$ into $\mathbb{H}^3/\rho(\pi _1(M))$.  Since $(i_{\rho})_*^{-1} \circ \rho: \pi _1(M) \longrightarrow \pi_1(M_{\rho})$ is an isomorphism  between aspherical manifolds there is a homotopy equivalence $h_{\rho}:M \longrightarrow M_{\rho}$ such that $(h_{\rho})_*$ is conjugate to $(i_{\rho})_*^{-1} \circ \rho$.  Then $\Theta$ maps $\rho$ to the marked homeomorphism class $[(M_{\rho}, h_{\rho})]$.  By the result of McCullough, Miller and Swarup above $\Theta$ is well-defined.  Marden's Isomorphism Theorem and Stability Theorem  ~\cite{Marden} and a result of Sullivan's \cite{Sullivan} implies that two convex co-compact representations are in the same component of $int(AH(\pi _1(M)))$ if and only if they have the same image under $\Theta$. 

\vspace{3mm}

{\bf The Klein-Maskit combination theorems}

\vspace{2mm}

Very roughly, the Klein-Maskit combination theorems give sufficient conditions for a group generated by two Kleinian groups to be Kleinian and tells us the topology of the manifold associated to the result.   The versions we state here are special cases of those given in \cite{Maskit}.

\vspace{4mm}

{\bf Theorem 1. (The First Klein-Maskit Combination Theorem)} {\it Let $J=<z \mapsto z+1>$ be a rank 1 parabolic subgroup of discrete groups $G_1$ and $G_2$ with $J\neq G_i$, $i=1,2$.  Assume that there is a horizontal line W which separates $\mathbb{C}$ into two closed disks $B_1$ and $B_2$ with  $\Lambda(G_j)\subset B_j$, and with the property that $B_j$ is precisely invariant under $J$ in $G_{3-j}$. Set $G=<G_1, G_2>$.   Then }

\vspace{2mm}
\renewcommand{\theenumi}{\arabic{enumi}}
\begin{enumerate}

\item {\it G is discrete;}\\

\item {\it $G=G_1 *_J G_2$;}\\

\item {\it  $G$ is geometrically finite if and only if both $G_1$ and $G_2$ are geometrically finite;}\\

\item {\it Let $C$ be the plane spanned by $W$ in $\mathbb{H}^3$ and let $B_i ^3$ be the half space bounded by $C$ and $B_i$ in $\mathbb{H}^3$, $i=1,2$.  Then $\mathbb{H}^3/G$ is isometric to the manifold obtained by  removing $B_1 ^3/J$ from $\mathbb{H}^3/G_2$ and removing $B_2 ^3/J$ from $\mathbb{H}^3/G_1$  and identifying the resulting manifolds along their common boundary $C/J$.}

\vspace{2mm}

\item {\it $({\mathbb H}^3 \cup \Omega(G))/G$ is orientation preserving homeomorphic to the manifold obtained by gluing $({\mathbb H}^3 \cup \Omega(G_1))/G_1$ to $({\mathbb H}^3 \cup \Omega(G_2))/G_2$ by identifying the punctured disk $B_1/J \subset \Omega (G_1)/G_1$ with the punctured disk $B_2/J \subset \Omega(G_2)/G_2$ via an orientation reversing homeomorphism.}\\

\end{enumerate}
\vspace{3mm}

The Second Combination Theorem deals with HNN-extensions of Kleinian groups.

\vspace{2mm}

{\bf Theorem 2. (The Second Klein-Maskit Combination Theorem)} {\it Let $J=<z\mapsto z+1>$ be a subgroup of a group $H$.   Suppose $f$ is a parabolic transformation $f(z)=z+c$ with $c$ not completely real.  Let  $A=\{z \ | \ a<\Im z< b\}$ with $|a-b|<|\Im c|$ and suppose that each component of $\mathbb{C}-A$ is precisely invariant under $J$ in $H$.  Set $G=<H, f>$.  Then}

\vspace{2mm}

\renewcommand{\theenumi}{\arabic{enumi}}
\begin{enumerate}

\item {\it $G$ is discrete};\\

\item {\it $G=H*_f$;}\\

\item {\it $G$ is geometrically finite if and only if $H$ is;}\\

\item {\it Let $B_1^3$ be the totally geodesic half-space whose closure meets $\widehat {\mathbb{C}}$ exactly in $B_1=\{ z \ | \Im z \geq b\}$, and $B_2 ^3$ the totally geodesic half-space whose closure meets $\widehat {\mathbb{C}}$ exactly in $B_2=\{ z \ | \Im z \leq a\}$.  Let $P_1$ be the plane in $\mathbb{H}^3$ $\{ (z,t) \ | \Im z=b\}$ and let $P_2$ be the plane  $\{ (z,t) \ | \Im z=a\}$.  Then $\mathbb{H}^3/G$ is isometric to the manifold obtained by removing $B_1^3/J$ and $B_2 ^3/J$ from $\mathbb{H}^3/H$ and identifying the resulting boundaries $P_1 /J$ and $P_2 /J$ via $f$.\\}  

\vspace{2mm}

\item {\it $({\mathbb H}^3 \cup \Omega(G))/G$ is orientation preserving homeomorphic to the manifold obtained from $({\mathbb H}^3 \cup \Omega(H))/H$ by identifying the punctured disk $B_1/J$ with the punctured disk $B_2/J$ in $\Omega(H)/H$ via the map induced by $f$.}\\

\end{enumerate}

\vspace{3mm}
\newpage
{\bf The Hyperbolic Dehn Surgery Theorem}

\vspace{3mm}

Let $\widehat M$ be a compact, irreducible, oriented 3-manifold whose boundary contains a single torus $T$.  There may be other surfaces in the boundary, but only one torus.  Choose a meridian $m$ and a longitude $l$ for the torus $T$ and consider them as a basis for $\pi _1(T)$.  For a pair of relatively prime integers $(p,q)$ denote by $\widehat {M}(p,q)$ the manifold obtained by performing $(p,q)$ Dehn surgery on $\widehat M$.  That is, $\widehat {M}(p,q)$ is obtained by sewing a solid torus $V$ to $\widehat M$ along $T$ via an orientation reversing homeomorphism which maps the meridian of $V$ to a simple closed curve in the homotopy class of $m^p l^q$ on $T$.  The following generalization of Thurston's Hyperbolic Dehn Surgery Theorem is due to T. Comar.  See also Bonahon and Otal \cite{Bonahon/Otal}.

\vspace{4mm}

{\bf Theorem 3. The Hyperbolic Dehn Surgery Theorem}(\cite{Comar}). {\it Let $\widehat M$ be a compact, oriented 3-manifold with one toroidal boundary component $T$.  Let $\widehat {N}=\mathbb{H}^3/\widehat {\Gamma}$ be a geometrically finite hyperbolic 3-manifold and $\phi:int(\widehat M)\longrightarrow \widehat N$ an orientation preserving homeomorphism between the interior of $M$ and $\widehat N$.  Further assume that every parabolic element of $\widehat \Gamma$ is conjugate to an element of $\phi _*(\pi _1(T))$.  Let $\{(p_n, q_n)\}$ be a sequence of distinct pairs of relatively prime integers.}

\vspace{2mm}

{\it Then, for all sufficiently large n, there exists a representation $\beta _n:\widehat {\Gamma} \longrightarrow PSL_2(\mathbb{C})$ with discrete image such that:
\renewcommand{\theenumi}{\arabic{enumi}}
\begin{enumerate}

\item $\beta _n (\widehat {\Gamma})$ is convex co-compact and uniformizes $\widehat {M}(p_n, q_n)$;\\

\item the kernel of $\beta _n$ is normally generated by $m^{p_n}l^{q_n}$; and\\

\item  $\{\beta _n\}$ converges to the identity representation of $\widehat \Gamma$.\\
\end{enumerate}
\vspace{2mm}

Moreover, if we let $i_n$ denote the inclusion of $\widehat M$ into $\widehat {M}(p_n, q_n)$, then there exists an orientation-preserving homeomorphism $\phi _n:int(\widehat {M}(p_n, q_n))\longrightarrow \mathbb{H}^3/\beta _n(\widehat {\Gamma})$ such that $\beta _n \circ \phi _*$ is conjugate to $(\phi _n)_* \circ (i_n)_*$.}

\vspace{5mm}

{\bf Some Notation:}  Let $A^{b}_{c}$ denote the horizontal strip $\{ z\in \mathbb{C}\ | \ c\leq \Im z \leq b\}$.

\vspace{2mm}

Let $H_{a}$ denote the half-plane $\{ z\in \mathbb{C}\ | \ \Im z \geq a\}$ and let $H^*_{a}$ denote the lower half-plane $\{z \in \mathbb{C}\ | \ \Im z \leq a\}$.

\vspace{2mm}

For $a \in \mathbb{C}$ let $\xi _a$ be the conformal transformation given by $\xi _a(z)=z+a$.

\vspace{2mm}

{\bf Lemma 1  }{\it Let $K\geq 1$.  Let $G$ be a Fuchsian group containing $\xi _1$ as a primitive element.  Then there is a constant $c=c(K,G)$ with the property that if $f$  is any $K$-quasiconformal deformation of $G$ fixing $0,1$ and $\infty$ then $H_c$ and $H^*_{-c}$ are precisely $J$-invariant in $f\circ G \circ f^{-1}$.}

\vspace{3mm}

{\bf Proof}

\vspace{2mm}

Suppose there does not exist a $c$ such that $H_c$ is precisely $J$-invariant in $f\circ G \circ f^{-1}$ for every  normalised $K$-quasiconformal deformation $f$ of $G$.  Then there is a sequence of normalized $K$-quasiconformal maps $\{f_n\}$ inducing quasiconformal deformations of $G$ and a sequence of group elements $\{g_n\}$ of $G-J$ such that $f_n\circ g_n \circ f_n ^{-1}(H_n)\cap H_n \neq \emptyset$.  Equivalently $g_n\circ f_n ^{-1}(H_n)\cap f_n ^{-1}(H_n)\neq \emptyset$.  Since $G$ is Fuchsian then Prop. VI.A.6 in \cite{Maskit} gives us that the half-spaces $H_1$ and $H^*_{-1}$ are precisely $J$-invariant in $G$.  Hence $f_n^{-1}(H_n)$ must intersect $A^1 _{-1}$.  

\vspace{2mm}

For each $n$, let $z_n$ be a point in $H_n$ such that $f_n ^{-1}(z_n)$ is in $A^1_{-1}$.  Since each $f_n$ commutes with $\xi _1$ we may assume that the $f_n ^{-1}(z_n)$ all lie in the (compact) rectangle $R=\{z \ | \ 0\leq \Re z \leq 1, \ -1\leq \Im z \leq 1\}$.  $\{f_n^{-1}(z_n)\}$ has an accumulation point $p\in R$, so we pass to a subsequence and reindex so that $f_n^{-1}(z_n)\longrightarrow p$. 

\vspace{2mm}

The space of all $K$-quasiconformal maps normalized to fix $0,1$ and $\infty$ is compact with respect to the topology of uniform convergence on compact subsets (\cite{Lehto} Theorems 2.1, 2.2).  Hence there is a normalized $K$-quasiconformal map $f$ such that, up to subsequence, $f_n\longrightarrow f$ uniformly on $R$.  Thus $f(p)=\lim _{n\rightarrow \infty}f_n(f_n ^{-1}(z_n))=lim _{n\rightarrow \infty} z_n = \infty$, supplying us with our contradiction.

\vspace{2mm}

By repeating the above argument or by appealing to symmetry we have that there is also a $c'$ such that $H^*_{-c'}$ is precisely $J$-invariant in $f\circ G \circ f^{-1}$ for any normalized $K$-quasiconformal deformation $f$ of $G$.  The desired constant $c(K,G)$ is the maximum of the two constants.$\Box$  

\vspace{6mm}

{\bf Corollary 2} {\it With the above assumptions $\Lambda(f \circ G \circ f^{-1})\subset A_{-c}^{c}\cup \{\infty\}$.}

\vspace{4mm}

{\bf Lemma 3  }{\it Let $G_1$,...,$G_N$ be quasifuchsian groups containing $J=<\xi_1>$ and suppose $\exists c>0$ such that $H^* _{-c}$ and $H_c$ are precisely $J$-invariant in each $G_i$.

\vspace{2mm}

Set $\Gamma$ to be the group generated by $\{\xi_{a_j i}G_j\xi_{a_j i}^{-1}\}^N _{j=1}$, where $0\leq a_1<a_2<\ldots <a_N$ are any real numbers with $a_{i+1}-a_{i}\geq 2c$ $\forall$ $i=1,\ldots, N-1$.

\vspace{2mm}

Then 

\renewcommand{\labelenumi}{\theenumi.}
\begin{enumerate}
\item $\Gamma = \xi_{a_1 i}G_1 \xi_{a_1 i}^{-1}*_J \xi_{a_2 i}G_2 \xi _{a_2 i}^{-1}*_J \cdots *_J \xi_{a_N i}G_N \xi _{a_N i}^{-1}$; \label{a}

\item $\Gamma$ is discrete and geometrically finite; \label{b}

\item $H_{a_N +c}\cup H^*_{a_1 -c}$ is precisely $J$-invariant in $\Gamma$; \label{c}

\item Let $S_j$ be a surface with boundary such that $G_j$ uniformizes $S_j \times I$ and for each $j$ let $\Delta _j$ be a component of $\partial S_j$.  Represent a solid torus $V$ as $D^2 \times S^1$ and form a manifold $M$ by attaching each $S_j \times I$ to $V$ by identifying $\Delta _j \times I$ with $[e^{2\pi i (4j-1)/4N}, e^{2\pi i (4j+1)/4N}] \times S^1$ via an orientation reversing homeomorphism.

\vspace{3mm}

Then $\Gamma$ uniformizes $M$, and if $\partial S_j=\Delta _j$ $\forall j$ then $({\mathbb H}^3 \cup \Omega(\Gamma))/\Gamma$ is orientation preserving homeomorphic to $M-\delta$, where $\delta$ is the simple closed curve $\{(1,0)\}\times S^1 \subset \partial M$ adjacent to $S_1$ and $S_N$.
\end{enumerate}

} 

\vspace{2mm}

Note: Since both $H_{a_N +c}$ and $H^*_{a_1 -c}$ are $J$-invariant, \ref{c}) implies that $H_{a_N +c}$ and $H^* _{a_1 -c}$ are precisely $J$-invariant in $\Gamma$.

\vspace{3mm}

{\bf Proof}

\vspace{2mm}

We proceed by induction on $N$.  Suppose that $N=1$, so that $\Gamma=\xi_{a_1 i}G_1 \xi _{a_1 i}^{-1}$.   $\Lambda(\Gamma)$ is a Jordan curve contained in $A^{a_1+c}_{a_1 -c}\cup \{\infty\}$ and the complement of this curve consists of two simply connected components each of which is stabilized by $\Gamma$.  Hence for $g\in \Gamma -J$, $g(H_{a_1+c})\subset H_{a_1 -c}$ and $g(H^*_{a_1 -c})\subset H^*_{a_1+c}$.  Since both $H_{a_1 +c}$ and $H^*_{a_1 -c}$ are precisely $J$-invariant in $\Gamma$ this implies that for $g\in \Gamma -J$, $g(H_{a_1+c})\subset H_{a_1 -c}\cap H^*_{a_1 +c}=A^{a_1 +c}_{a_1 -c}$ and  $g(H^*_{a_1 -c})\subset H^*_{a_1+c}\cap H_{a_1 -c}=A^{a_1 +c}_{a_1 -c}$ and hence that $g(H_{a_1 +c}\cup H^*_{a_1 -c})\cap (H_{a_1 +c}\cup H^*_{a_1 -c})=\emptyset$ proving that part 3) holds and hence, since 1), 2) and 4) are clearly true, the base step of our induction.

\vspace{2mm}

Suppose now that $G=<\xi_{a_1 i}G_1\xi_{a_1 i}^{-1},...,\xi _{a_{k-1}i}G_{k-1}\xi _{a_{k-1}i}^{-1}>$ and $H=\xi_{a_k i}G_k \xi _{a_k i}^{-1}$.  We prove that the lemma is true for $\Gamma=<G,H>$.

\vspace{2mm}

To prove parts 1) and 2) consider the line $W=\xi _{(a_{k-1}+c)i}({\mathbb R})$.  $W$ divides $\widehat {\mathbb C}$ into two disks $H_{a_{k-1}+c}$ and $H^*_{a_{k-1}+c}$.  By our inductive hypotheses $H_{a_{k-1}+c}$ is precisely $J$-invariant in $G$ and by assumption $H^*_{a_{k-1}+c}$ is precisely $J$-invariant in $H$ since $H^*_{a_k-c}\supset H^*_{a_{k-1}+c}$ and $H^*_{a_k -c}$ is precisely $J$-invariant in $H$.  Parts 1) and 2) now follow with an application of Theorem 1.

\vspace{2mm}

To prove part 3) first observe the following:

\vspace{2mm}

\hspace{5mm}for $g\in G-J$, $g(H_{a_k +c}\cup H^*_{a_1 -c})\subset g(H_{a_{k-1}+c}\cup H^*_{a_1 -c})\subset A^{a_{k-1}+c}_{a_1 -c}$ by our inductive assumption;

\vspace{2mm}

\hspace{5mm}for $h\in H-J$, $h(H_{a_k +c}\cup H^*_{a_1 -c})\subset h(H_{a_k +c}\cup H^*_{a_k -c})\subset A^{a_k +c}_{a_k -c}$;

\vspace{2mm}

\hspace{5mm}for $g\in G-J$, $g(A^{a_k +c}_{a_k -c})=g(H_{a_k -c}\cap H^*_{a_k +c})\subset g(H_{a_{k-1}+c})\subset A^{a_{k-1} +c}_{a_1 -c}$;and

\vspace{2mm}

\hspace{5mm}for $h\in H-J$, $h(A^{a_{k-1}+c}_{a_1 -c})\subset h(H^*_{a_k -c})\subset A^{a_k +c}_{a_k -c}$.

\vspace{2mm}

Let ${\cal T}_G$ be a complete set of right coset representatives for $J$ in $G$ containing the identity  and similarly let ${\cal T}_H$ be a complete set of right coset representatives for $J$ in $H$ containing the identity.  Then any element $g$ of $\Gamma$ can be expressed as $g=g_1 h_1\cdots g_n h_n\xi ^t$ where $g_i \in {\cal T}_G$, $h_i \in {\cal T}_H$ and $t\in {\mathbb Z}$.  Using the above observations it is clear that $g(H_{a_k +c} \cup H^*_{a_1 -c})\subset A^{a_k +c}_{a_1 -c}$.

\vspace{2mm}

Hence it follows that for $g\in \Gamma -J$, $g(H_{a_k +c}\cup H^*_{a_1 -c})\cap (H_{a_k +c}\cup H^*_{a_1 -c})=\emptyset$ and part \ref{c}) of the lemma follows.

\vspace{2mm}

It only remains to prove part 4).

\vspace{2mm}

Let $C_i$ be the totally geodesic plane meeting $\widehat{\mathbb C}$ in $\{ z \ | \ \Im z=a_i +c\}$.  Let $B^- _{i}$ be the half-space bounded by $C_i$ and $H_{a_{i-1} +c}$ and let $B^+ _{i}$ be the half-space bounded by $C_i$ and $H^*_{a_{i-1}+c}$.  Then by Theorem 1 applied inductively,

\begin{equation*}
\begin{split}
 {\mathbb H}^3/\Gamma &\cong ({\mathbb H}^3 -B^-_1)/\xi_{a_1 i}G_1\xi_{a_1 i}^{-1}\cup  ({\mathbb H}^3-(B^+_1 \cup B^-_2))/\xi_{a_2 i}G_2\xi_{a_2 i}^{-1}\cup
\cdots\\ & \cdots \cup ({\mathbb H}^3-(B^+_{N-1} \cup B^-_N))/\xi_{a_{N-1} i}G_{N-1}\xi_{a_{N-1} i}^{-1}\cup ({\mathbb H}^3 -B^{+}_N)/\xi_{a_N i}G_N\xi_{a_N i}^{-1}
\end{split}
\end{equation*}
\vspace{2mm}
For each $i$ let $A_i$ be an annulus in $S_i$ with $\Delta _i$ as one of its boundary components.  Then by Theorem 1, part 5) the right hand side of the above expression is the interior of
\begin{equation*}
\begin{split}
M'=(S_1 \times I)\cup (S_2 \times I)\cup \cdots \cup (S_N \times I)
\end{split}
\end{equation*}

\vspace{2mm}

\noindent where in the union $A_i \times \{1\}$ is identified with $A_{i+1}\times \{0\}$, $i=1,\ldots, N-1$ by  orientation reversing homeomorphisms.

\vspace{2mm}

We claim that $M'$ is orientation preserving homeomorphic to $M$.  Consider the solid torus $\bigcup A_i \times I$, where again the appropriate identifications are made as above. The removal of this solid torus from the above manifold is orientation preserving homeomorphic to the disjoint union $\bigsqcup _{i=1}^N (S_i \times I - \Delta_i \times I)$.  Hence $M'$ is obtained by attaching $\bigsqcup _{i=1}^N S_i \times I$ to a solid torus along annuli $\Delta _i \times I$.  In other words, $M'$ is orientation preserving homeomorphic to $M$.

\vspace{2mm}

A similar reasoning shows that $({\mathbb H}^3 \cup \Omega(\Gamma))/\Gamma$ is of the desired form.
$\Box$
\vspace{5mm}

{\bf III: The Proof of Theorem A}

\vspace{2mm}

Our goal is to prove the following result:

\vspace{2mm}

{\bf Theorem A} {\it There exists a geometrically finite representation $\rho \in AH(\pi _1(M_k))$ with $\Theta(\rho)=[(M_k, id)]$ such that for each $[(M_k^{\tau}, h_{\tau})]\in {\cal A}(M)$ there is a sequence of convex co-compact representations $\rho _n ^{\tau}$ converging to $\rho$ with $\Theta(\rho_n^{\tau})=[(M_k ^{\tau},h_{\tau})]$.}
\vspace{3mm}

{\bf Remark}  It is clear that Theorem A follows from Theorem B.  For expository reasons we prove Theorem A first and then generalize the result in Theorem B.  Theorem B states that if you choose a quasiconformal constant $K\geq 1$ then there is a geometrically finite representation such that its $K$-quasiconformal deformations are contained in the closure of every component of $CC(\pi _1(M_k))$.  Theorem A furnishes us with this representation.  So while there is no mention of a $K$ in the statement of Theorem A the proof makes use of the $K$ of the statement of Theorem B.
\vspace{2mm}

{\bf Proof}

\vspace{2mm}

 The idea is that we find a ``nice'' uniformization for $M_k$, call it $\Gamma _k$,  and ``shuffling parabolics'' $g_{\tau}$ so that the group $\widehat \Gamma _k ^{\tau}$ generated by $\Gamma _k$ and $g_{\tau}$  is a geometrically finite uniformization of $\widehat M_k ^{\tau}$. For each $\tau$ we have an inclusion  $j_{\tau}: \Gamma _k \longrightarrow \widehat {\Gamma _k ^{\tau}}$.  Using The Hyperbolic Dehn Surgery Theorem we obtain convex co-compact representations  $\beta _n ^{\tau}:\widehat {\Gamma _k ^{\tau}}\longrightarrow PSL_2(\mathbb{C})$ with image uniformizing $M_k ^{\tau}$ in the correct homeomorphism class and converging to the identity representation.  Suppose $\phi : int(M_k) \longrightarrow \mathbb{H}^3/\Gamma _k$ is an orientation preserving homeomorphism.  We will consider the representations $\rho _n ^{\tau}:\pi _1(M_k) \longrightarrow PSL_2(\mathbb{C})$ given by 

\vspace{2mm}

\[ \rho _n ^{\tau}=\beta _n ^{\tau}\circ j_{\tau} \circ \phi _*\]

\vspace{2mm}

$\rho _n ^{\tau}$ converges to $\rho_K= j_{\tau}\circ \phi _*=\phi_*$ with image $\Gamma _k$, a geometrically finite uniformization of $M_k$.

\vspace{2mm}

{\bf Step 1:}{\it The uniformization of $M_k$ and $\widehat {M_k ^{\tau}}$}

\vspace{2mm}

For each $j \in 1,2,...,k$ let $G_0 (j)$ be a Fuchsian model for $B(j)$ containing $\xi _1$ as a primitive element. Fix $K\geq 1$ and let $c_i=c(K, G_0 (j))$ be the constant given by lemma 1, $i=1,..,k$.  Set $C=max\{c_1,...,c_k\}$.  Choose, for each $\tau$ a prime $p_{\tau}$ so that $p_{\tau}=p_{\tau '}$ if and only if $\tau = \tau '$ and so that $p_{\tau}>4Ck$, $\forall \tau$.

\vspace{2mm}

By appealing to the Chinese Remainder Theorem or simply by constructing them, we find integers $d_{\tau}$, one for each $\tau$, with the following properties:

\vspace{2mm}
\hspace{35mm}$d_{\tau}\equiv 1 \mod{p_{\tau}}$ and

\hspace{35mm}$d_{\tau}\equiv 0 \mod{p_{\tau '}}$ for $\tau \neq \tau '$.

Define integers $a_j$ by

\vspace{2mm}

\[ a_j=\sum _{\tau} d_{\tau}(jp_{\tau} + 4C\tau ^{-1}(j))\]

\vspace{2mm}

Set $G (j)=\xi _{a_j i}G_0 (j) \xi _{a_j i}^{-1}$ and define $\Gamma _k$ to be the group generated by $G (1),.., G(k)$.  By Lemma 3 $\Gamma _k$ is discrete, geometrically finite and uniformizes $M_k$.

\vspace{2mm}

Set $g_{\tau}=\xi _{p_{\tau}i}$ and $\widehat {\Gamma _k ^{\tau}}=<\Gamma _k, g_{\tau}>$.  We wish to show that $\widehat {\Gamma _k ^{\tau}}$ is a geometrically finite uniformization of $\widehat {M_k ^{\tau}}$.

\vspace{2mm}

Observe now that $a_j\equiv 4C\tau ^{-1}(j) \mod{p_{\tau}}$ and so that $G(j)$ is conjugate by powers of $g_{\tau}$ in $\widehat {\Gamma _k ^{\tau}}$ to $\xi _{4C\tau ^{-1}(j) i}G_0(j)\xi _{4C\tau ^{-1}(j)i}^{-1}$  The latter is equal to $\xi _{4Cli}G_0(\tau(l))\xi _{4Cli}^{-1}$ for some $l$.  Hence $g_{\tau}$ has ``shuffled'' the limit sets of the $G(j)$'s in the manner that $\tau$ shuffles $1,2,...,k$ (see figure 2).  Denote the group generated by $\xi _{4Ci} G_0(\tau(1))\xi _{4Ci}^{-1}$, $\xi _{8Ci} G_0(\tau(2))\xi _{8Ci}^{-1}$,..., $\xi _{4Cki} G_0(\tau(k))\xi _{4Cki}^{-1}$ by $\Gamma _k ^{\tau}$.  Then $\widehat {\Gamma _k ^{\tau}}$ is generated by $\Gamma _k ^{\tau}$ and $g_{\tau}$.  Applying Lemma 3 we find that $\Gamma _k ^{\tau}$ is geometrically finite and uniformizes $M_k ^{\tau}$.

\vspace{5mm}

\fbox{\rule[-4cm]{0cm}{4cm}\epsfig{file=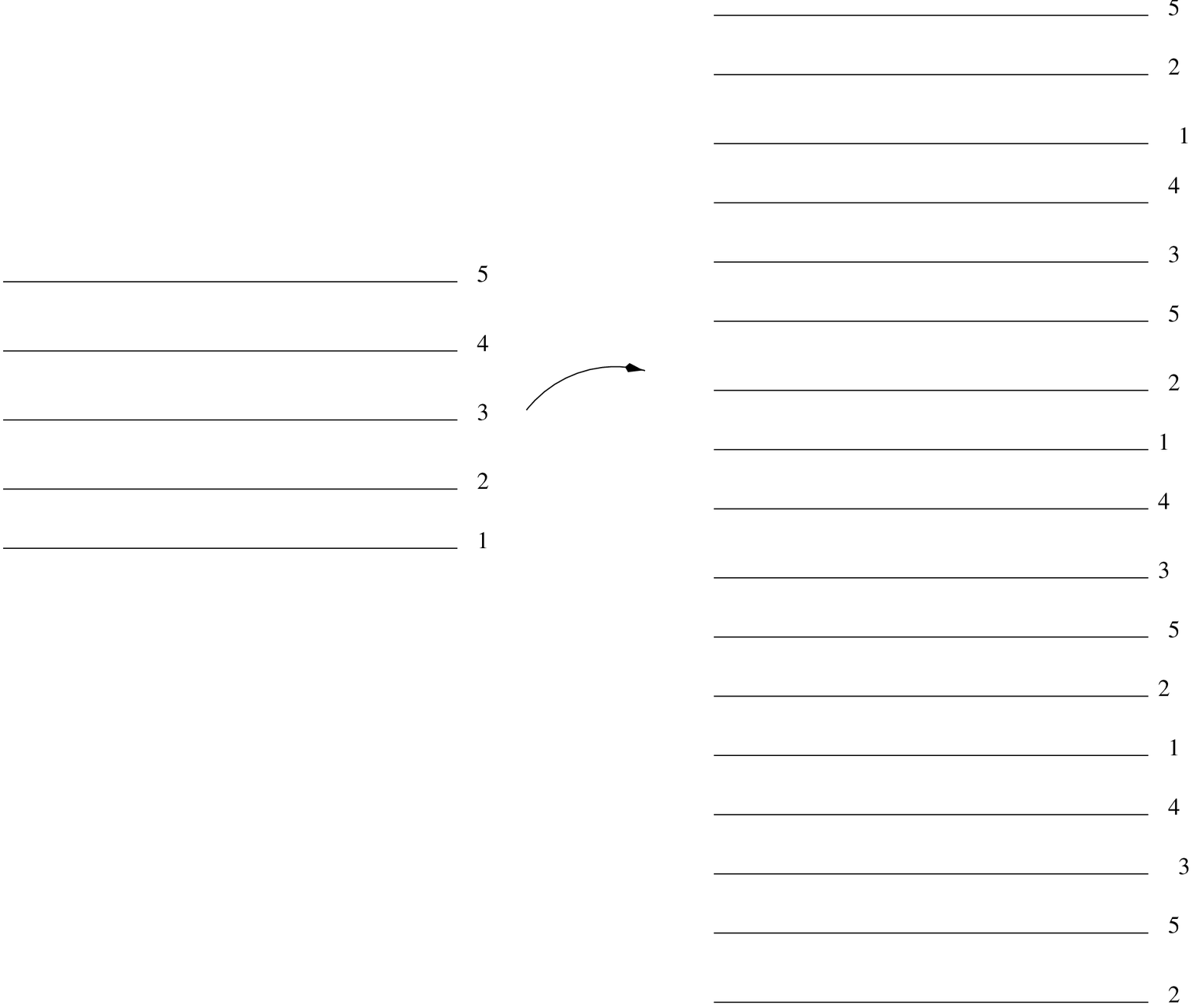, height=10cm, width = 7cm, angle=-90}}

\vspace{2mm}
\begin{center}
Figure 2\\
The shuffling depicted here is for $\tau = (13)(24)$ with $k=5$.\\
\end{center}

Lemma 3 gives us that both $H_{4Ck+3C}$ and $H^* _{3C}$ are precisely $J$-invariant in $\Gamma _k ^{\tau}$ and $\Lambda (\Gamma _k ^{\tau})\subset A^{4Ck +C}_{3C}\subset A^{4Ck +3C}_{3C}$.  Since $p_{\tau}>4Ck$ we can apply Theorem 2 to conclude that ${\widehat {\Gamma _k ^{\tau}}}$ is discrete, torsion free and geometrically finite. 

We shall see that $\widehat {\Gamma _k ^{\tau}}$ uniformizes ${\widehat {M_k ^{\tau}}}$. $({\mathbb H}^3\cup \Omega(\Gamma _k ^{\tau}))/\Gamma _k ^{\tau}$ is orientation preserving homeomorphic to $M_k ^{\tau} - \delta$ where $\delta \subset \partial V$ is a longitudinal simple closed curve adjacent to both $B(\tau (1))$ and $B(\tau(k))$.  Take two distinct annuli in $(\partial M_k ^{\tau}-\delta)\cap \partial V$ each having $\delta$ as a boundary component.  By Theorem 2, $({\mathbb H}^3\cup \Omega(\widehat {\Gamma _k ^{\tau}}))/{\widehat {\Gamma _k ^{\tau}}}$ is orientation preserving homeomorphic to the manifold obtained from $M_k ^{\tau} - \delta$ by identifying these two annuli.  This resulting manifold is clearly orientation preserving homeomorphic to $M_k ^{\tau} - \alpha$, where $\alpha$ is a core curve for $V$.  Hence

\begin{center}
${\mathbb H}^3/\widehat {\Gamma _k ^{\tau}} \cong int ({\mathbb H}^3 \cup \Omega(\widehat {\Gamma _k ^{\tau}})/{\widehat {\Gamma _k ^{\tau}}} \cong M_k ^{\tau} - {\cal N}(\alpha) \cong \widehat {M_k ^{\tau}}$\hspace{8mm}
\end{center}

\vspace{5mm}
\begin{center}
\fbox{\rule[-6cm]{0cm}{6cm}\epsfig{file=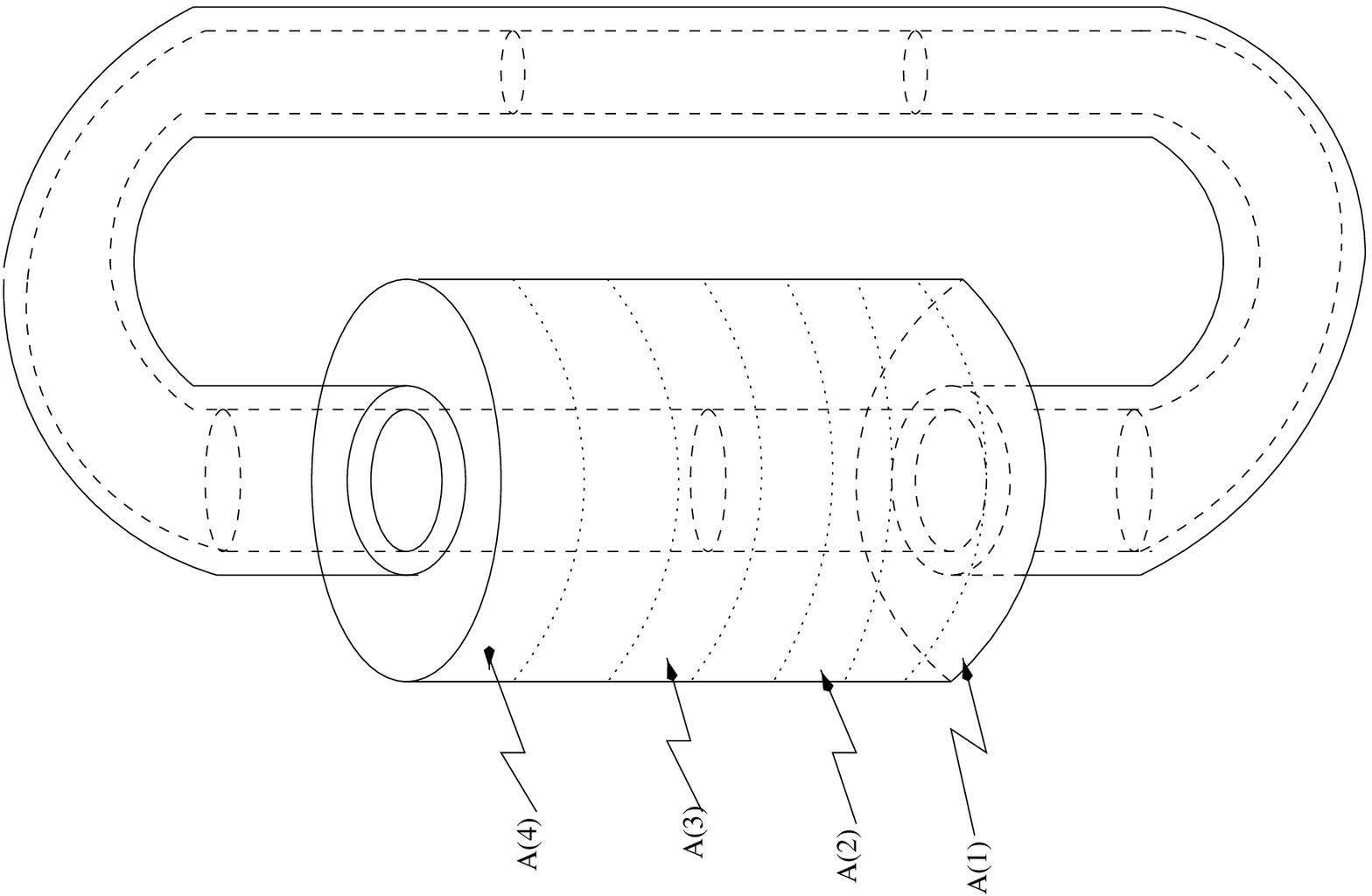, height=8cm, width=10cm, angle=-90}}\\
\vspace{2mm}
Figure 3.\\ 
\end{center}
\vspace{3mm}

{\bf Step 2:}{\it Constructing the algebraic limits}

\vspace{2mm}
Let $\alpha$ be the core curve of $V$ and consider $M_k ^{\tau} - {\cal N}(\alpha)$, which is homeomorphic to ${\widehat {M_k ^{\tau}}}$.  Choose a meridian $m$ and a longitude $l$ on the resulting torus boundary component so that $l$ is homotopic to $\alpha$ in $V$ and $m$ bounds a disk in ${\cal N}(\alpha)$.  By Lemma 10.3 of \cite{Anderson/Canary/McCullough} the manifold ${\widehat {M_k ^{\tau}}}(1,n)$ obtained by performing $(1,n)$ Dehn surgery  on ${\widehat {M_k ^{\tau}}}$ is homeomorphic to $M_k ^{\tau}$, $\forall n$.  Moreover, if $q_n:\widehat{M_k ^{\tau}}(1,n)\longrightarrow M_k ^{\tau}$ is the homeomorphism furnished by Lemma 10.3 of \cite{Anderson/Canary/McCullough} then $q_n$ can be chosen to be the identity on each $B(j)\subset \widehat{M_k ^{\tau}}(1,n)$.

\vspace{2mm}

Choose an orientation preserving homeomorphism $\phi ^{\tau}:int(\widehat M_k ^{\tau})\longrightarrow \mathbb{H}^3/\widehat \Gamma _k ^{\tau}$.  The Hyperbolic Dehn Surgery Theorem supplies us with representations $\beta _n ^{\tau}:\widehat \Gamma _k ^{\tau}\longrightarrow PSL_2(\mathbb{C})$ such that for large enough $n$ 

\vspace{2mm}

\renewcommand{\theenumi}{\roman{enumi}}
\begin{enumerate}
\item $\beta _n ^{\tau}(\widehat \Gamma _k ^{\tau})$ is convex co-compact; 
\item there exist orientation preserving homeomorphisms $\phi_n ^{\tau}:int(\widehat M_k ^{\tau})\longrightarrow \mathbb{H}^3/\beta _n ^{\tau}(\widehat \Gamma _n ^{\tau})$ such that if $i_n:\widehat M_k ^{\tau}\longrightarrow \widehat M_k ^{\tau}(1,n)$ is inclusion then $\beta _n ^{\tau}\circ (\phi  ^{\tau})_*$ is conjugate to $(\phi _n ^{\tau})_* \circ (i_n)_*$.
\end{enumerate}

\vspace{2mm}

Choose an orientation preserving homeomorphism $\phi: int(M_k) \longrightarrow \mathbb{H}^3/\Gamma _k$.  Let $j_{\tau}:\Gamma _k \longrightarrow \widehat \Gamma _k ^{\tau}$ be inclusion and set

\vspace{2mm}
\[ \rho _n ^{\tau}=\beta _n ^{\tau} \circ j_{\tau} \circ \phi _*\]
\vspace{2mm}

The kernel of $\beta _n ^{\tau}$ is normally generated by $g_{\tau}\xi _1 ^n$.  Hence it is the group generated by the set of all elements of the form $hg_{\tau}\xi _1 ^n h^{-1}$ with $h\in \widehat {\Gamma _k ^{\tau}}$.  $\Gamma _k$ has trivial intersection with this group, and hence $\rho _n ^{\tau}$ is faithful.

\vspace{2mm}

For ease of notation denote $\xi _{4C\tau ^{-1}(j) i}G (j)\xi _{4C\tau ^{-1}(j)i}^{-1}$ by $G^{\tau}(j)$ , $j=1,...,k$.  $j_{\tau}(\xi _{a_j i}G(j) \xi _{a_j i}^{-1})= g_{\tau}^{n_j}G^{\tau}(j) g_{\tau}^{-n_j}$ where $n_j = \frac {a_j - 4C\tau^{-1}(j)}{p_{\tau}}$.  Since $\beta _n ^{\tau}(g_{\tau}^{n_j}G^{\tau}(j)g_{\tau}^{-n_j})=\beta_n ^{\tau}(G^{\tau}(j))$ and the image of these $G^{\tau}(j)$'s generates $\beta _n ^{\tau}(\widehat {\Gamma _k ^{\tau}})$ we have that the image of $\rho _n ^{\tau}$ is that of $\beta _n ^{\tau}$ and hence is discrete, convex co-compact, and uniformizes $M_k ^{\tau}$.

\vspace{2mm}

Recall the map $\Theta: AH(\pi_1 (M_k))\longrightarrow {\cal A}(M_k)$ defined previously.  

\vspace{2mm}

We wish to show that $\Theta (\rho _n ^{\tau})=[(M_k ^{\tau}, h_{\tau})]$ where $h_{\tau}$ is a homotopy equivalence between $M_k$ and $M_k ^{\tau}$ which is the identity off of $V$.  This will imply that $\rho _n ^{\tau}$ belongs to the component of $CC(\pi _1(M_k))$ corresponding to $[(M_k ^{\tau}, h_{\tau})]$.  Let $q_n:{\widehat {M_k ^{\tau}}}(1,n)\longrightarrow M_k ^{\tau}$ be an orientation preserving homeomorphism and $p_{\tau}:\mathbb{H}^3/\Gamma _k \longrightarrow \mathbb{H}^3/\widehat \Gamma _k ^{\tau}$ an orientation preserving  covering map with $(p_{\tau})_* =j_{\tau}$.  We claim that $r_{\tau}=q_n\circ i_n \circ (\phi ^{\tau})^{-1}\circ p_{\tau}\circ \phi$ is a homotopy equivalence between $int(M_k)$ and $int(M_k ^{\tau})$.  Since $\rho _n ^{\tau}$ is conjugate to $(\phi _n ^{\tau}\circ (q_n)^{-1})_*\circ (r_{\tau})_*$ with $\rho _n ^{\tau}$ faithful and both $\phi _n ^{\tau}$ and $q_n$ homeomorphisms, $(r_{\tau})_*$ is injective. By similar reasoning $(r_{\tau})_*$ is surjective.  

\vspace{2mm}
Denote the homotopy inverse of $h_{\tau}$ by $\overline{h_{\tau}}$.  The composition $h=r_{\tau}\circ \overline{h_{\tau}}$ is a homotopy equivalence from $M_k ^{\tau}$ to itself.  If $K(M)$ denotes the characteristic submanifold of $M$ (see \cite{Jaco} or \cite{Johannson} for details) then by Theorem 24.2 of \cite{Johannson} $h$ is homotopic to a map (again called $h$) such that $h:K(M_k^{\tau})\longrightarrow K(M_k^{\tau})$ is  homotopy equivalence and $h:\overline{M_k^{\tau}-K(M_k^{\tau})}\longrightarrow \overline{M_k^{\tau}-K(M_k ^{\tau})}$ is a homeomorphism.  The characteristic submanifold of $M_k ^{\tau}$ is $M_k ^{\tau}-\bigcup{\cal N}(A(i))$ and so consists of the $B(j)$ together with a solid torus $V_0\subset V$ such that $\partial V_0 \cap \partial V\subset \partial M_k^{\tau}\cap \partial V$ (see Section 4 of \cite{Culler/Shalen}).  Applying Proposition 28.4 of \cite{Johannson} we can homotope $h$ to a map (again, called $h$) that is a homeomorphism when restricted to each $B(j)$ and to $V_0$ by a homotopy that preserves the $\partial _0(B(j))$ and the lids of $B(j)$, for $j=1,\ldots,k$.  There is the possibility that this new map reverses the orientation on $\partial _0(B(j))$, for each $j$.  If this were the case then $h_*(\alpha)$ would not be conjugate to $\alpha$ in $\pi _1(M_k ^{\tau})$, where $\alpha$ generates $\pi _1(V)$.  However, $(r_{\tau})_*\circ ({\overline {h_{\tau}}})_*(\alpha)$ is conjugate to $\alpha$ and hence it is not the case that $h$ reverses the orientation on $\partial _0(B(j))$. Hence $h$ may be assumed to be the identity on $\partial _0(B(j))$ $\forall j$.  Since $h$ is homotopic to the identity on $V_0$, $h$ is an orientation preserving homeomorphism from $M_k ^{\tau}$ to itself.  In particular $[(M_k^{\tau},r_{\tau})]=[(M_k^{\tau},h_{\tau})]$.

\vspace{2mm}

Observe that since $i:M_k ^{\tau}\longrightarrow \widehat M_k ^{\tau}$ is an embedding and $i_n \circ i$ is homotopic to an orientation preserving homeomorphism $i_n\circ i (M_k^{\tau})$ is a compact core for $\widehat M_k ^{\tau}(1,n)$.  Homotoping $i_n \circ i$ slightly we may assume that the image of $M_k ^{\tau}$ lies in the interior of $\widehat M_k ^{\tau} (1,n)$.  Then $\phi _n ^{\tau}\circ i_n \circ i(M_k ^{\tau})$ is a compact core for $\mathbb{H}^3/\rho _n ^{\tau}(\pi _1 (M_k))$.  Hence 

\begin{align*}
\Theta(\rho _n ^{\tau}) &= [(\phi _n ^{\tau}\circ i_n \circ i (M_k ^{\tau}), \phi _n ^{\tau}\circ q_n ^{-1}\circ r_{\tau})]\\
&= [(i_n \circ i(M_k ^{\tau}), q_n ^{-1}\circ r_{\tau})]\\
&= [(q_n \circ i_n \circ i(M_k ^{\tau}), r_{\tau})]\\
&= [(M_k ^{\tau}, r_{\tau})]\\
&= [(M_k ^{\tau}, h_{\tau})]\\
\end{align*}

and this second to last equality holds because $q_n \circ i_n \circ i$ is homotopic to an orientation preserving homeomorphism.

\vspace{2mm}

Hence $\rho _n ^{\tau}$ belongs to the component of $CC(\pi _1(M_k))$ indexed by $[(M_k ^{\tau}, h_{\tau})]$ as claimed.

\vspace{2mm}

Finally, $\rho_n ^{\tau}$ converges to $\rho _K = \phi_*$ which can be seen to  belong to the closure of the component of $CC(\pi _1(M_k))$ indexed by $[(M_k, id)]$ by either setting $\tau =id$ and using the proof above, or by a direct application of the Hyperbolic Dehn Surgery Theorem (see Remark(1) following the proof of Theorem B in \cite{Anderson/Canary/McCullough}), or by appealing to corollary 6 of Ohshika (\cite{Ohshika}).   

\vspace{2mm}

We have completed the proof.  $\Box$

\vspace{5mm}

{\bf Theorem B}(Intersection Contains A connected, Uncountable Set)

\vspace{2mm}

 {\it For every $K\geq 1$ there is a geometrically finite representation $\rho _K$ of $\pi _1(M_k)$ such that the set $O_K\subset QC(\rho _K)$ consisting of all of the $K$-quasiconformal deformations of $\rho _K$ is contained in the closure of every component of $CC(\pi _1 (M_k))$ and $\Theta(\rho_K)=[(M_k, id)]$.}

\vspace{2mm}

{\bf Proof}

\vspace{2mm}

Fix $\rho '$ in $O_K$ and suppose that $f$ is the $K$-quasiconformal map inducing the isomorphism between $\rho '$ and $\rho = \rho _K$, where $\rho _K$ is the representation constructed in Theorem A.  By conjugating if necessary we assume that $f$ fixes $0,1$ and $\infty$, so that $J$ is a primitive subgroup of $\rho '(\pi _1(M_k))$.

\vspace{2mm}

Set $G(j) '=\rho '(\rho ^{-1}(G(j)))$.  Then $G(j)'$ is quasi-Fuchsian for each $j$ and, by Corollary 2, $\Lambda(G(j)')\subset A^{a_j +C}_{a_j -C}$ (these $a_j$ are the $a_j$ that appeared in the construction of $\Gamma _k$).  Let $\Gamma _k '= \rho ' (\pi _1(M_k))$ and for each $\tau \in S_k/\mathbb{Z}_k$ set ${\widehat \Gamma _{\tau}}'=<\Gamma _k ', g_{\tau}>$.  We will show that ${\widehat \Gamma _{\tau}}'$ is discrete, torsion-free, geometrically finite and that  the quotient of $\mathbb{H}^3$ by ${\widehat \Gamma _{\tau}}'$ is  homeomorphic (via an orientation preserving homeomorphism) to the interior of $\widehat M_k ^{\tau}$.

\vspace{2mm}

As in Theorem A, define $n_j = \frac {a_j - 4C\tau ^{-1}(j)}{p_{\tau}}$ and set $G^{\tau}(j)'=g_{\tau}^{-n_j}G (j)'g_{\tau}^{n_j}$.  Set $\Gamma _{\tau}'$ to be the group generated by the $G^{\tau}(j)'$.  Then ${\widehat \Gamma _{\tau}}'$ is generated by $\Gamma _{\tau}'$ and $g_{\tau}$.  By Lemma 3 $\Gamma _{\tau}'$ uniformizes $M_k ^{\tau}$.  Just as in Theorem A we can apply Lemma 1, Lemma 3 and Theorem 2 to show that $\widehat \Gamma _{\tau} '$ is geometrically finite and uniformizes $\widehat M_k ^{\tau}$.

\vspace{2mm}

We now proceed exactly as in Theorem A.  To begin, appealing to Theorem 3 (The Hyperbolic Dehn Surgery Theorem) gives a sequence of representations $\alpha _n ^{\tau}$ of $\widehat \Gamma _{\tau}'$ in $PSL_2(\mathbb{C})$ converging to the identity, so that the image of each $\alpha _n ^{\tau}$ is convex co-compact uniformization of $M_k ^{\tau}$.  Then, with $j_{\tau}':\Gamma _k ' \longrightarrow {\widehat \Gamma _{\tau}}'$ being inclusion, we set

\vspace{2mm}

\[ \rho ^{\tau '}_{n}  =\alpha _n ^{\tau}\circ j_{\tau}' \circ \rho '\]

\vspace{2mm}

Then, as in Theorem A, for $n$ sufficiently large  each $\rho ^{\tau '}_{n}$ is faithful and has discrete image uniformizing $M_k ^{\tau}$ convex co-compactly.  Moreover, $\rho ^{\tau '}_{n}$ converges to $\rho '$.$\Box$

\vspace{5mm}
Observing that we could have easily done the above proof substituting $M_k^{\sigma}$ whenever we saw $M_k=M_k^{id}$, for some $\sigma$ in $S_k$, we have the following corollary to the proof.

\vspace{2mm}

{\bf Corollary C}(There are many sets in the intersection)

\vspace{2mm}

{\it For every $K\geq 1$ and $[(M_k^{\sigma}, h_{\sigma})]\in {\cal A}(M_k)$ there is a geometrically finite representation $\rho _K ^{\sigma}$ and a set $O_K ^{\sigma}\subset QC(\rho _K ^{\sigma})$ consisting of the $K$-quasiconformal deformations of $\rho _K ^{\sigma}$ such that $O_K ^{\sigma}$ is contained in the closure of every component of $CC(\pi _1 (M_k))$ and $\Theta (\rho _K ^{\sigma})=[(M_k ^{\sigma}, h_{\sigma})]$.}

\vspace{2mm}
This naturally leads us to ask if the intersection of the closures of the path components of $CC(\pi _1(M_k))$ is disconnected.

\vspace{2mm}
\bibliographystyle{plain}
\bibliography{bump}
\end{document}